\documentclass[12pt]{article}
\usepackage{amsfonts,amsmath,amsxtra}
\usepackage{amssymb}
\usepackage{bm}
\usepackage[utf8]{inputenc}
\usepackage{lscape}
\usepackage[normalem]{ulem}
\usepackage[russian,english]{babel}
\usepackage{xcolor}
\usepackage{pst-circ}
\usepackage{tikz}
\numberwithin{equation}{section}

\def\a{\alpha}

\def\C{\mathbb{C}}

\def\bl{{\bm{\lambda}}}
\def\bL{{\bm{\Lambda}_n}}
\def\beq{\begin{equation}}
	\def\eeq{\end{equation}}
\def\beqq{\begin{equation*}}
	\def\eeqq{\end{equation*}}
	\def\bn{{\bm{\nu}}}
\def\bs{\begin{split}}
	\def\es{\end{split}}

\def\bl{\boldsymbol{\lambda}}
\def\bu{\boldsymbol{u}}
\def\bv{\boldsymbol{v}}

\def\d{\partial}

\def\g{\gamma}

\def\HS{{H}}
\def\i{\hat{i}}

\def\l{\lambda}
\def\la{\lambda}

\def\nn{[n]}
\def\nnn{[n-1]}

\def\R{\mathbb{R}}
\def\Re{\mathrm{Re}\,}
\def\Res{\operatorname{Res}}

\def\tS{\tilde{S}}

\def\vf{\varphi}

\def\x{\bm{x}}
\def\z{\bm{z}_n}

\newtheorem{lemma}{Lemma}
\newtheorem{proposition}{Proposition}

\newtheorem{theorem}{Theorem}
\newtheorem{remark}{Remark}
\newcommand{\rf}[1]{(\ref{#1})}

\begin{document}
\begin{center}
\phantom{1}

\vspace{-2cm}
\hfill ITEP-TH-21/21\\
\hfill IITP-TH-16/21\\ [20mm]
	{\bf \large Wave function for  $GL(n,\R)$ hyperbolic Sutherland model II.

\smallskip
{  \,    Dual Hamiltonians }}
	\bigskip

{\bf S. Kharchev$^{\,\star\,\natural}$,\, S. Khoroshkin$^{\,\star\,\circ, \natural}$,
}\medskip\\
$^\star${\it Institute for Theoretical and Experimental Physics, Moscow, Russia;}\smallskip\\
$^\natural${\it
	Institute for Information Transmission Problems RAS (Kharkevich Institute),Bolshoy Karetny per. 19, Moscow, 127994, Russia;}\smallskip\\
$^\circ${\it National Research University Higher School of Economics, Moscow, Russia.}
\end{center}
\begin{abstract}\noindent
Recently we found Mellin-Barnes integrals, representing the wave
 function for $GL(n,\R)$ hyperbolic Sutherland model. In present paper,
 we establish bispectral properties of this wave function with
 respect to dual Ruijesenaars-Macdonald operators.
\end{abstract}	

	\bigskip
\section{Introduction}

 {\bf 1.	}
 In \cite{KK} we found Barnes type integrals presenting the wave function  for  $GL(n,\R)$ hyperbolic Sutherland Hamiltonian
 with arbitrary coupling constant $g>0$
 \beq\label{i1}
H^{(g)}_{2}(x_1,\ldots,x_n)=-\sum_{i=1}^n\frac{\d^2}{\d x_i^2}+\sum_{{ i\not=j}}\frac{g(g-1)}{\sh^2(x_i-x_j)}.
\eeq	
It has a form
\beq\label{i2}
\Psi^{(g)}_{\l_1,\ldots, \l_n}(x_1,\ldots,x_n)
=\prod_{j<k}\sh^{g}|x_j-x_k|\times \Phi^{(g)}_{\l_1,\ldots, \l_n}(x_1,\ldots,x_n),
\eeq
where
\beq\label{i3}
\begin{split} \Phi^{(g)}_{\l_1,\ldots, \l_n}&(x_1,\ldots,x_n)
=\\
 \int\limits_{\imath\R^{\frac{n(n-1)}{2}}}\prod\limits_{i=1}^{n-1}&
 \frac{\prod\limits_{j=1}^{i}\prod\limits_{k=1}^{i+1}
 \Gamma\left(\frac{\g_{i,j}-\g_{i+1,k}+g}{2}\right) \Gamma\left(\frac{\g_{i+1,k}-\g_{i,j}+g}{2}\right)}
 {\prod\limits_{1\leq r\neq s\leq i}\Gamma\left(\frac{\g_{i,r}-\g_{i,s}}{2}\right)
 \Gamma\left(\frac{\g_{i,r}-\g_{i,s}+2g}{2}\right)}
	 	e^{\sum_{i,j=1}^n(\g_{i,j}-\g_{i-1,j})x_i}
	 	\prod\limits_{\stackrel{i=1}{j\leq i}}^{n-1}d\g_{i,j}.
	 \end{split}
	 \eeq	
	 Here $\l_i=\g_{n,i}$, $\l_i\in \imath\R$, $g>0$, and $\g_{i,j}=0$ if $i<j$.
The function 	 $\Phi^{(g)}_{\l_1,\ldots, \l_n}(x_1,\ldots,x_n)$ is analytical in a strip around $\R^n\subset\C^n$ and presumably
coincides with the corresponding Heckman--Opdam hypergeometric series, invariant with respect to permutations of the coordinates $x_k$. Relations  \rf{i3} supplement integral formulas by M. Halln\"{a}s and S. Ruijesenaars \cite{HR1}, where the  wave functions  for  $GL(n,\R)$ hyperbolic Sutherland system are given by integrals over space variables.

 {\bf 2. } On the other hand, it is known that spectral problems for integrable systems with Hamiltonians expressed via
 Sekiguchi-Macdonald operators \cite{M,S} admit bispectral extensions. In the paper \cite{Ch} O. Chalykh posed a bispectral problem for the hyperbolic
 Sutherland problem with integer coupling constant and found its solution in a form of Baker--Akhiezer function (see \cite{CE} { for the further development}).
 Later M. Noumi, J. Shiraishi \cite{NS} set a bispectral problem for Ruijsenaars--Macdonald $q$--difference operators. They proved the existence of  a
 solution of this problem in a form of basic hypergeometric series in a region $|x_1|\gg|x_2|\gg\cdots\gg|x_n|$ and studied analytical problems
  of this solution. Besides, they observed  certain invariance properties of their solution under the change $t\leftrightarrow qt^{-1}$ of Macdonald parameters.
 { In \cite{LNS} these results were generalized to elliptic  Ruijsenaars difference operators}.

 {\bf 3. } Set
{ \beq\label{i4}
 {\cal H}_r^{(g)}(\l_1,\ldots,\l_n)=(-1)^{r(n-1)}
 \sum_{I_r\subset \nn}\,\prod_{i\in I_r,j\in\nn\setminus I_r}
 \frac{\l_{i}-\l_{j}+2-2g}{\l_{i}-\l_{j}}\cdot
\prod_{i\in I_r} T_{\l_{i}},
 \eeq
 where the sum is taken over all subsets $I_r$ of the set $\nn=\{1,\ldots,n\}$ of cardinality $r$; and  $T_{\l_k} f(\l_k)=f(\l_k+2)$.} Commuting operators \rf{i4}
  are known as Ruijsenaars--Macdonald difference operators \cite{M,R}, { describing Ruijsenaars rational relativistic generalization of Calogero-Sutherland model}.
  Our main result is
  \begin{theorem} The function \rf{i2} solves the spectral problem for Ruijsenaars--Macdonald difference operators \rf{i4} with arbitrary coupling constant $g>1$, namely,
  	\beq\label{i5}
   {\cal H}_r^{(g)}(\l_1,\ldots,\l_n)\Psi^{(g)}_{\l_1,\ldots, \l_n}(x_1,\ldots,x_n)=e_r\left(e^{2x_1},\ldots, e^{2x_n}\right)\Psi^{(g)}_{\l_1,\ldots, \l_n}(x_1,\ldots,x_n).
   \eeq	
  	\end{theorem}
  Here $e_r(y_1,\ldots,y_n)$, $r=1,\ldots,n$ is $r$-th elementary symmetric function,
  \beq \label{i5a}e_r(y_1,\ldots,y_n)=\sum_{i_1<\ldots< i_r}y_{i_1}\cdots y_{i_r}.\eeq
   The wave function \rf{i2}, originally defined for imaginary values of parameters $\l_i$, is analytical in a strip around real subspace $\R^n\subset\C^n$ and then analytically continued to other values of  parameters $\l_i$. On the other hand, shifts in difference operators \rf{i4} are performed in real direction, which means that the relation \rf{i5} is understood in a sense of analytical continuation or as  an integral equation with Cauchy type kernels.

  Besides, we have by \cite{KK} the relations
  \beq\label{i6}\begin{split}
  	H_1^{(g)}(x_1,\ldots,x_n)	\Psi_{\l_1,\ldots,\l_n}^{(g)}(x_1,\ldots,x_n)=&\Big(\sum_p \l_p\Big)\Psi_{\l_1,\ldots,\l_n}^{(g)}(x_1,\ldots,x_n),\\
  	H_2^{(g)}(x_1,\ldots,x_n)	\Psi_{\l_1,\ldots,\l_n}^{(g)}(x_1,\ldots,x_n)=&-\Big(\sum_p \l_p^2\Big)\Psi_{\l_1,\ldots,\l_n}^{(g)}(x_1,\ldots,x_n),
  \end{split}\eeq	
 Here $	H_2^{(g)}(x_1,\ldots,x_n)$ is Sutherland Hamiltonian \rf{i1} and $	H_1^{(g)}(x_1,\ldots,x_n)$ is proportional to the total impulse operator, $H_1^{(g)}(x_1,\ldots,x_n)=\frac{\d}{\d x_1}+\ldots + \frac{\d}{\d x_n}$.
 Thus the function \rf{i2} solves the bispectral problem for operators $H_r^{(g)}(x_1,\ldots,x_n)$, $r=1,2$, and
 ${\cal H}_r^{(g)}(\l_1,\ldots,\l_n)$, $r=1,\ldots,n$.

  Note that  M. Halln\"{a}s and S. Ruijesenaars also obtained integral presentations for the wave functions of relativistic Ruijsenaars  system \cite{HR2, HR3}. So far we did not establish precise connection of the rational degeneration of these presentations with the formula \rf{i3}.

   The plan of the paper is as follows. In Section 2 we recall basic formulas of the theory of Macdonald operators,
    describe certain gauge equivalent Macdonald operators, and introduce the notions of dual Macdonald operators and bispectral problem for their wave functions. Besides, we recall difference equations on the weights of measures, associated to Macdonald operators, which guarantee symmetricity of these operators. In Section 3 we specialize these notions to hyperbolic Sutherland model, performing appropriate limit. We describe corresponding dual Hamiltonians (closely related to Ruijsenaars integrable system), and measures, related to them. Note that the measures we use do not coincide with the formal limit of measures, used for Macdonald polynomials. The corresponding weights differ by periodic functions and could be rather  regarded as deformations of Sklyanin measure \cite{GKL}. Then we prove that the wave function constructed in \cite{KK} solves the spectral problem for dual Hamiltonians. The proof is a  play with contours in Mellin-Barnes integrals with the use of related identities on rational functions. Note that here the restriction $g>1$ appears, contrary to \cite{KK}, where the coupling constant is arbitrary positive.

    Finally, in Appendix we prove certain identities on rational functions which serve the proof of our main result.

 \section{Macdonald operators}
{\bf 1. Measures and gauges. }
  Recall some basic facts about Macdonald operators. For $q,t\in\C^\times$ and ${\z}=(z_1,\ldots,z_n)$ the Macdonald operator
  $ M_r=M_r(\z|q,t)$, $r=1,\ldots,n$ is the operator in the space of analytical in $\left(\C^\times\right)^n$ functions on $n$ variables, given by the
  expressions
{ \beq\label{m0}
 M_r(\z|q,t)=\sum_{I_r\subset \nn}\,\prod_{i\in I_r,j\in\nn\setminus I_r}
 \frac{tz_{i}-t^{-1}z_j}{z_{i}-z_j}\cdot\prod_{i\in I_r} T_{z_{i}},
 \eeq
 where the sum is taken over all subsets $I_r$ of the set $\nn=\{1,\ldots,n\}$ of cardinality $r$} and  $T_{z_i}$ is a shift operator
 \beqq
 (T_{z_i}f)(z_1,\ldots,z_i,\ldots,z_n)=f(z_1,\ldots,q^2z_i,\ldots,z_n).
 \eeqq
 These operators commute \cite{M}:
 $$[M_r,M_s]=0,$$
 are invariant with respect to permutations of the coordinates, and Macdonald polynomials form
 an eigenbasis for the action of these operators in the space of symmetric polynomials.

 \bigskip\noindent
{ Assume that $|q|<1$ and introduce the standard notations $(z;q^2)_\infty=\prod_{i=0}^\infty(1-zq^{2i})$. Consider the function
 \beq\label{aa1}
 \phi_{a,b}(\z):=\prod_{j\neq k}
\frac{(az_j/z_k;q^2)_\infty}{(bz_j/z_k;q^2)_\infty}
 \eeq
where $a,b\in\C^\times$.  It satisfies the  relations
 \beq\label{aa2}
T_{z_i}\circ\phi_{a,b}(\z)=
\prod_{j\neq i}\frac{(z_i-aq^{-2}z_j)(bz_i-z_j)}{(az_i-z_j)(z_i-bq^{-2}z_j)}\,
\phi_{a,b}(\z)\circ T_{z_i}.
 \eeq
We use this function for two important particular choices of parameters $a$ and $b$.
\medskip

 Set $a=1$, $b=t^2$. The function
$\phi_{1,t^2}(z)$ coincides with the weight  $\Delta_{q,t}(\z)$
 \beq\label{m1}
 \Delta_{q,t}(\z)=\prod_{j\neq k}
 \frac{(z_j/z_k;q^2)_\infty}{(t^2z_j/z_k;q^2)_\infty},
 \eeq
 for the scalar product in the space of analytical functions.
 (named "another scalar product" in \cite{M}):
 \beq \label{m3}
 (f,g)_{q,t}=\oint_{|z|_i=1}f(\z)\overline{g(\z)} \Delta_{q,t}(\z)\prod_i\frac{dz_i}{z_i}
\eeq
It is known \cite{M} that for  real  $q$ and $t$ Macdonald operators \rf{m0} are symmetric with respect to the scalar product \rf{m3}.
 This follows from the relation
 \beq \label{m4}
 T_{z_i}\circ\Delta_{q,t}(\z)=\prod_{j\neq i}\frac{qz_i-q^{-1}z_j}{z_i-z_j}
 \cdot\frac{tz_i-t^{-1}z_j}{qt^{-1}z_i-q^{-1}tz_j}\,
 \Delta_{q,t}(\z)\circ T_{z_i}.
\eeq

Put  $a=q,\,b=q^2t^{-2}$ and consider the function \cite{NS}
\beq\label{m5}
	\phi(\z)=\prod_{j\neq k}
	\frac{(qz_j/z_k;q^2)_\infty}{(q^{2}t^{-2}z_j/z_k;q^2)_\infty}
\eeq
\begin{lemma}\label{lemma0}
	Operators $M_r(\z|q,t)$ and $M_r(\z|q,qt^{-1})$ are gauge equivalent:
	\beq\label{m6}
		M_r(\z|q,t)\circ\phi(\z)=\phi(\z)\circ M_r(\z|q,qt^{-1}).\phantom{\int}\hspace{-0.5cm}
	\eeq
\end{lemma}
 {\it Proof}. Direct calculation { using the commutation relations
\beq\label{aa4}
T_{z_i}\circ\phi(\z)=
\prod_{j\neq i}\frac{qt^{-1}z_i-q^{-1}tz_j}{tz_i-t^{-1}z_j}\,
\phi(\z)\circ T_{z_i},
\eeq
which are particular cases of (\ref{aa2})} \hfill{$\square$}

\smallskip\noindent
 {\bf 2.  Duality}.
Consider the spectral problem
 \beq\label{m9}
 M_r(\z|q,t)\Phi_\bL(\z)=e_r(\bL)\Phi_\bL(\z),
 \eeq
 where $\bL=(\Lambda_1,\ldots,\Lambda_n)$ is a tuple of complex parameters, $e_r(\bL),\,(i=1,\ldots,n)$ are elementary symmetric functions \rf{i5a}. It is known from the theory of quantum Knizhnik-Zamolodchikov equation
 \cite{C,MS,TV} that this problem may be extended to a bispectral problem on meromorphic functions on complex variables $\z$ and $\bL$. Moreover M. Noumi and J. Shiraishi \cite{NS} constructed by iterated summations of basic hypergeometric series meromorphic solutions of the bispectral problem
 \beq\label{m10}
\begin{split}
	&M_r(\z|q,t)\Phi_\bL(\z)=e_r(\bL)\Phi_\bL(\z),\\
 	&M_r(\bL|q,t)\Phi_\bL(\z)=e_r(\z)\Phi_\bL(\z),
 \end{split}
\eeq
 where operators $M_r(\bL|q,t)$ have the same form as initial operators
 $M_r(\z|q,t)$
{ \beq\label{m11}
 M_r(\bL|q,t)=\sum_{I_r\subset \nn}\,\prod_{i\in I_r,j\in\nn\setminus I_r}
 \frac{t\Lambda_{i}-t^{-1}\Lambda_{j}}
 {\Lambda_{k}-\Lambda_j}\cdot
\prod_{i\in I_r} T_{\Lambda_{i}},
 \eeq }
 where difference operators $T_{\Lambda_i}$ are given by the relation:
 $$T_{\Lambda_i}f(\Lambda_1,\ldots,\Lambda_i,\ldots,\Lambda_n)=
 f(\Lambda_1,\ldots,q^2\Lambda_i,\ldots,\Lambda_n).$$
 Due to equivalence \rf{m6} one may consider another bispectral system
  \beq\label{m10a}
 \begin{split}
 	M_r(\z|q,t)\Phi_\bL(\z)=&e_r(\bL)\Phi_\bL(\z),\\
 	M_r(\bL|q,qt^{-1})\Phi_\bL(\z)=&e_r(\z)\Phi_\bL(\z),
 \end{split}
 \eeq
 with
  \beq\label{m11a}
 M_r(\bL|q,qt^{-1})=\sum_{I_r\subset\nn}\,\prod_{i\in I_r,j\in\nn\setminus I_r}
 \frac{qt^{-1}\Lambda_{i}-q^{-1}t\Lambda_{j}}
 {\Lambda_{i}-\Lambda_j}\cdot\prod_{i\in I_r}T_{\Lambda_{i}}.
 \eeq
 We call operators $M_r(\bL|q,t)$ and gauge equivalent operators $M_r(\bL|q,qt^{-1})$ {\it dual }
Macdonald operators.
 In analogy with \rf{m6}, they are connected by the relation
		\beq\label{aa7}
		M_r(\bL|q,t)\circ\phi(\bL)=\phi(\bL)\circ M_r(\bL|q,qt^{-1}).\phantom{\int}\hspace{-0.5cm}
		\eeq
		where the function
		\beq\label{aa8}
		\phi(\bL)=\prod_{j\neq k}
		\frac{(q\Lambda_j/\Lambda_k;q^2)_\infty}
		{(q^{2}t^{-2}\Lambda_j/\Lambda_k;q^2)_\infty}
		\eeq
		satisfies the commutation relations
		\beq\label{aa9}
		T_{\Lambda_i}\circ\phi(\bL)=
		\prod_{j\neq i}\frac{qt^{-1}\Lambda_i-q^{-1}t\Lambda_j}
		{t\Lambda_i-t^{-1}\Lambda_j}\,
		\phi(\bL)\circ T_{\Lambda_i}
		\eeq
 The corresponding weight functions
  \beq\label{m7}
 \Delta_{q,t}(\bL)=\prod_{j\neq k}
 \frac{(\l_j/\l_k;q^2)_\infty}{(t^2\l_j/\l_k;q^2)_\infty},\qquad \text{and} \qquad
 \Delta_{q,qt^{-1}}(\bL)=\prod_{k\neq j}
 \frac{(\l_j/\l_k;q^2)_\infty}{(q^2t^{-2}\l_j/\l_k;q^2)_\infty}
 \eeq satisfy the difference equations
\beq \label{mn4}
\begin{split}
& T_{\Lambda_i}\circ\Delta_{q,t}(\bL)=
 \prod_{j\neq i}\frac{q\Lambda_i-q^{-1}\Lambda_j}{\Lambda_i-\Lambda_j}
 \cdot\frac{t\Lambda_i-t^{-1}\Lambda_j}{qt^{-1}\Lambda_i-q^{-1}t\Lambda_j}\,
 \Delta_{q,t}(\bL)\circ T_{\Lambda_i},\\
& T_{\Lambda_i}\circ\Delta_{q,qt^{-1}}(\bL)=
 \prod_{j\neq i}\frac{q\Lambda_i-q^{-1}\Lambda_j}{\Lambda_i-\Lambda_j}
 \cdot\frac{qt^{-1}\Lambda_i-q^{-1}t\Lambda_j}{t\Lambda_i-t^{-1}\Lambda_j}\,
 \Delta_{q,qt^{-1}}(\bL)\circ T_{\Lambda_i}
 \end{split}
\eeq

\smallskip\noindent
{\bf 3. Limit to Sutherland} { Let $q=e^{\pi\imath\tau},\,{\rm Im}\,\tau>0$, $t=q^g,\ g>0$,
 and tend $\tau$ to $0$.} Then the expansion of $M_1(z_1,\ldots,z_n|q,t)$ with respect to $\tau$ looks as follows
 \beq\label{m12}
 M_1(\z|q,t)=n+\pi\imath\tau\HS_1^{(g)}(\z)-\frac{(\pi\tau)^2}{2}
 \HS_2^{(g)}(\z)+o(\tau^2),
 \eeq
 where
 \beq \label{m13}\begin{split}
 	\HS_1^{(g)}(\z)=&2\sum_{i=1}^nz_i\d_{ z_i},\\
 	\HS_2^{(g)}(\z)=&4\left(\sum_{i=1}^n\left(z_i\d_{ z_i}\right)+g\sum_{i<j}\frac{z_i+z_j}{z_i-z_j}\,(z_i\d_{z_i}-z_j\d_{z_j})+
 	\frac{1}{12}g^2n(n^2-1)\right)
 	\end{split}
 \eeq
 After the change of variables
 \beq\label{m14}z_i=e^{2x_i}\eeq
 formulas \rf{m13} transform to { reduced} Hamiltonians of hyperbolic Sutherland model:
  \beq \label{m15}\begin{split}
 	\HS_1^{(g)}(\x_n)=&\sum_{i=1}^n\d_{ x_i},\\
 	\HS_2^{(g)}(\x_n)=&\sum_{i=1}^n\d^2_{ x_i}+2g\sum_{i<j}\cth(x_i-x_j)\,(\d_{x_i}-\d_{x_j})+
 	\frac{1}{3}g^2n(n^2-1).
 \end{split}
 \eeq
 The difference equations \rf{m5} on the weight $\Delta_{q,t}(\z)$  of the scalar product turns in this limit
 to differential equations on the limiting weight $\Delta_g(\z)$:
 \beq \label{m16}
 z_i\d_{z_i}\Delta_g(\z)=g\sum_{j\not=i}\frac{z_i+z_j}{z_i-z_j}\Delta_g(\z)
 \eeq
 which solution is
 \beq\label{m17}
 \Delta_g(\z)=\prod_{j\not=i}\left|\sqrt{\frac{z_i}{z_j}}-\sqrt{\frac{z_j}{z_i}}\right|^{2g}
 \eeq
 so that in variables $\x$, see \rf{m14}
 \beq\label{m17a}
 \Delta_g(\x_n)=\prod_{j\not=i}\sh^{2g}|x_i-x_j|
 \eeq
 which coincides with the direct limit of the measure weight function for integer $g$.
 \section{{ Dual Hamiltonians }}

{\bf 1. Sutherland specialization of dual Macdonald operators}. Set
\beq\label{s1}
 \Lambda_i:=e^{\pi\imath\tau \la_i},\ \ \ {{\rm Im}\,\tau>0}.
\eeq
Then in the limit $\tau\to 0$  the dual Macdonald operators \rf{m11}  { and \rf{m11a}} tend to dual Sutherland Hamiltonians
 \beq\label{s2}
D_r^{(g)}(\bl_n)=
\sum_{I_r\subset \nn}\,\prod_{i\in I_r,j\in\nn\setminus I_r}
\frac{\l_{i}-\l_{j}+2g}{\l_{i}-\l_{j}}\cdot
\prod_{i\in I_r}T_{\l_{i}}
\eeq
{ and
 \beq\label{s2a}
D_r^{(1-g)}(\bl_n)=
\sum_{I_r\subset \nn}\,\prod_{i\in I_r,j\in\nn\setminus I_r}
\frac{\l_{i}-\l_{j}+2-2g}{\l_{i}-\l_{j}}\cdot
\prod_{i\in I_r}T_{\l_{i}}
\eeq }
where { $T_{\Lambda_k}:=T_{\l_k}=\exp\{2\d_{\l_k}\}$}.

\noindent
{ Consider the function
\beq
\vf(\bl_n)=\prod_{p\neq q}\Gamma\Big(\frac{\l_p-\l_q+2-2g}{2}\Big).
\eeq
{ This is the solution of equations
\beq
T_{\l_i}\circ\vf(\bl_n)=
(-1)^{n-1}\prod_{j\neq i}\frac{\l_i-\l_j-2g+2}{\l_i-\l_j+2g}\,
\vf(\bl_n)\circ T_{\l_i}
\eeq
which are the liming cases of relations (\ref{aa9}). }

\noindent
In analogy with (\ref{m6}) there is the gauge transformation
\beq\label{ns2}
D_r^{(g)}(\bl_n)\circ\vf(\bl)=\vf(\bl)\circ{\cal H}_r^{(g)}(\bl_n)
\eeq
where
\beq\label{s3}
\begin{split}
& {\cal H}_r^{(g)}(\bl_n)=(-1)^{r(n-1)}D_r^{(1-g)}(\bl_n)=\\
&=
(-1)^{r(n-1)}\sum_{I_r\subset \nn}\,\prod_{i\in I_r,j\in\nn\setminus I_r}
\frac{\l_{i_k}-\l_{j}+2-2g}{\l_{i_k}-\l_{j}}\cdot\prod_{i\in I_r}
T_{\l_i},
\end{split}
\eeq
}
It is known that these operators are gauge equivalent to   Hamiltonians of the rational degeneration of Ruijsenaars system constructed in \cite{R}.

In the limit $\tau\to 0$, the difference equations (\ref{mn4})
%on the limiting measure weights
become
\begin{align}
\label{s4}
&T_{\la_i}\circ\mu^{(g)}(\bl_n)=
\prod_{j\neq i}\frac{\la_i-\la_j+2}{\la_i-\la_j}\cdot
\frac{\la_i-\la_j+2 g}{\la_i-\la_j+2-2g}\,
\mu^{(g)}(\bl_n)\circ T_{\la_i},\\
&T_{\la_i}\circ\mu^{(1-g)}(\bl_n)=
\prod_{j\neq i}\frac{\la_i-\la_j+2}{\la_i-\la_j}\cdot
\frac{\la_i-\la_j+2-2 g}{\la_i-\la_j+2g}\,
\mu^{(1-g)}(\bl_n)\circ T_{\la_i}\label{s5}.
\end{align}
These equations admit the solutions
\begin{align}\label{s6}
&\mu^{(g)}(\bl_n)=\prod\limits_{j\not=k}
\Gamma^{-1}\Big(\frac{\lambda_j-\lambda_k}{2}\Big)
\Gamma^{-1}\Big(\frac{\lambda_j-\lambda_k}{2}+{1-g}\Big),\\
&\mu^{(1-g)}(\bl_n)=\prod\limits_{j\not=k}
\Gamma^{-1}\Big(\frac{\lambda_j-\lambda_k}{2}\Big)
\Gamma^{-1}\Big(\frac{\lambda_j-\lambda_k}{2}+{g}\Big). \label{s7}
\end{align}
 We implicitly use the solution \rf{s7} for dual Hamiltonians \rf{s3} and denote it by ${\mu}(\bl_n)$.
\begin{remark}
Note that measure weights \rf{s6} and \rf{s7} differ by periodic functions from formal limits of the measure weights
  \rf{m7}.
 \end{remark}
\begin{remark}
	At the point $g=1/2$ when the wave function of the hyperbolic Sutherland model can be given by $GL(n,\R)$ zonal spherical functions {\rm\cite{GKL}}, both measures are proportional to  Sklyanin measure
	\beq  \label{s8}{\mu}_S(\bl_n)=\prod\limits_{j\not=k}
	\Gamma^{-1}({\lambda_j-\lambda_k}{}) \eeq
	due to Legendre duplication formula {\rm\cite{BE}}
$$
\Gamma(2z)=\frac{2^{2z-1}}{\sqrt\pi} \Gamma(z)\Gamma\Big(z+\frac{1}{2}\Big).
$$
%\sout{\blue{ used for Toda chain systems.}}
	\end{remark}
Our main result, Theorem 1, states that the  { wave} function \rf{i3} satisfies the system of spectral equations with respect to Hamiltonians \rf{i4} and thus satisfies the system of bispectral equations \rf{i5}, \rf{i6}

\smallskip\noindent
{\bf 2. Proof of Theorem 1.} The basic formula \rf{i3} for the wave function of the reduced Sutherland Hamiltonian can be written in inductive way
\beq\label{p1}
\Phi_{\bl_n}^{(g)}({\bm x}_n)=
\int\limits_{\i\R^{n-1}}\mu(\bn_{n-1})d\bn_{n-1}K(\bl_n,\bn_{n-1})
e^{h_n(\bl_n,\bn_{n-1})x_n}
\Phi_{\bn_{n-1}}^{(g)}({\bm x}_{n-1}),
\eeq
where
\begin{align*}&\bl_n=\left(\l_{1},\ldots, \l_{n}\right), &&\bn_{n-1}=\left(\nu_{1},\ldots, \nu_{n-1}\right),\\
&\x_n=\left(x_1,\ldots,x_n\right), && \x_{n-1}=\left(x_1,\ldots,x_{n-1}\right),
\end{align*}
\beqq
h_n(\bl_n,\bn_{n-1})=\sum_{j=1}^n\l_{j}-\sum_{j=1}^{n-1}\nu_{j},
\eeqq
\beq\label{p2}
K(\bl_n,\bn_{n-1})=\prod\limits_{j=1}^{n-1}\prod\limits_{k=1}^{n}
\Gamma\left(\frac{\nu_{j}-\l_{k}+g}{2}\right) \Gamma\left(\frac{\l_{k}-\nu_{j}+g}{2}\right),
\eeq
and
 $\mu(\bn_{n-1})d\bn_{n-1}$ is modified Sklyanin measure \rf{s7}. Here
\beq\label{p3}
\mu(\bn_{n-1})=\prod\limits_{1\leq r\neq s\leq n-1}
\Gamma^{-1}\left(\frac{\nu_{r}-\nu_{s}}{2}\right)
\Gamma^{-1}\left(\frac{\nu_{r}-\nu_{s}+2g}{2}\right).
\eeq
The important property of the integration contour in the integral \rf{p1}:
 it separates two series of poles:
 \beq \label{p1a}\begin{split} \nu_{i}&=\l_{j}+g+2k,\qquad k=0,1,\ldots \qquad\text{and}
 	\\ \nu_{i}&=\l_{j}-g-2k,\qquad k=0,1,\ldots
\end{split} \eeq
\begin{proposition}\label{proposition1}
	For $g>1$ and any $r=1,\ldots,n$ we have the following relation
	\beq \label{p4}\begin{split}
	&{\cal H}_r^{(g)}(\bl_n)\Phi_{\bl_n}^{(g)}(\x_n)=\\ &e^{2x_N}\int\limits_{\i\R^{n-1}}\mu(\bn_{n-1})d\bn_{n-1}K(\bl_n,\bn_{n-1})
e^{h_n(\bl_n,\bn_{n-1})x_n}
{\cal H}_{r-1}^{(g)}(\bn_{n-1})\Phi_{\bn_{n-1}}^{(g)}({\bm x}_{n-1})+\\ \ \ &\int\limits_{\i\R^{n-1}}\mu(\bn_{n-1})d\bn_{n-1}K(\bl_n,\bn_{n-1})
e^{h_n(\bl_n,\bn_{n-1})x_n}
{\cal H}_r^{(g)}(\bn_{n-1})\Phi_{\bn_{n-1}}^{(g)}({\bm x}_{n-1}).
\end{split}	\eeq
\end{proposition}

\smallskip\noindent
 {\it Proof of Theorem 1} now consists of the
 (double) inductive use of  Proposition \ref{proposition1}. Indeed, for $n=1$ the wave
 function is an exponent
  $$
  \Phi_{\l_1}^{(g)}(x_1)=e^{\l_1x_1},$$
  while the Macdonald operator ${\cal H}_1^{(g)}(\l_1)=T_{\l_1}$. Evidently,
  $$T_{\l_1}\Phi_{\l_1}^{(g)}(x_1)=e^{2x_1}\Phi_{\l_1}^{(g)}(x_1).$$
  Then, using successively \rf{p4} for $r=1$, we get the equality
  \beq \label{p5}
  {\cal H}_1^{(g)}(\bl_n)\Phi_{\bl_n}^{(g)}(\x_n)=\left(e^{2x_1}+\ldots + e^{2x_n}\right)\Phi_{\bl_n}^{(g)}(\x_n)
  \eeq
for any $n\geq 1$. Then we proceed in the same manner for ${\cal H}_2^{(g)}(\bl_n)$, etc., ${\cal H}_n^{(g)}(\bl_n)$, using the inductive relations for
   elementary symmetric functions,
   $$e_r(z_1,\ldots, z_n)=z_ne_{r-1}(z_1,\ldots,z_{n-1})+e_{r}(z_1,\ldots,z_{n-1}).$$ \hfill{$\Box$}

   \smallskip\noindent
    {\it Proof of  Proposition \ref{proposition1}}.
    The left hand side of \rf{p4} looks as
      \beq\label{p6}\begin{split}
     &{\cal H}_r^{(g)}(\bl_n)\Phi_{\bl_n}^{(g)}(\x_n)=\\&(-1)^{r(n-1)}
    \sum_{I_r\subset \nn}\,\prod_{i\in I_r,j\in\nn\setminus I_r}
\frac{\l_{i_k}-\l_{j}+2-2g}{\l_{i_k}-\l_{j}}\cdot\prod_{i\in I_r}
T_{\l_i}\Phi_{\bl_n}^{(g)}(\x_n).
     \end{split}\eeq
Consider the summand corresponding subset
$I_r=\{i_1,i_2,\ldots, i_r\}$. Denote this summand by $J_{I_r}$:
{\beq \label{p7}\begin{split} J_{I_r}=
&\prod_{i\in I_r,j\in\nn\setminus I_r}
\frac{\l_{i_k}-\l_{j}+2-2g}{\l_{i_k}-\l_{j}}\cdot\prod_{i\in I_r}
T_{\l_i}\times \\
&\int\limits_{\i\R^{n-1}}\mu(\bn_{n-1})d\bn_{n-1}K(\bl_n,\bn_{n-1})
e^{h_n(\bl_n,\bn_{n-1})x_n}
\Phi_{\bn_{n-1}}^{(g)}({\bm x}_{n-1}),
\end{split}\eeq}
 Shifts do note touch the integration measure $\mu(\bn_{n-1})d\bn_{n-1}$, but act nontrivially on the kernel $K(\bl_n,\bn_{n-1})$,
  exponent $e^{h_n(\bl_n,\bn_{n-1})x_n}$ and on the integration contour. We have by \rf{p2}
{ \beq\label{p8}\begin{split}
%T_{\l_{i}}K(\bl_n,\bn_{n-1})=&\prod_{a=1}^{N-1}\frac{\l_i-\nu_a+g}{\nu_a-\l_i+g-2}K(\bl_n,\bn_{n-1}) ,\\
 \prod_{i\in I_r}
 T_{\l_{i}}K(\bl_n,\bn_{n-1})=&\prod_{i\in I_r}\prod_{a=1}^{N-1}\frac{\l_{i}-\nu_a+g}{\nu_a-\l_{i}+g-2}K(\bl_n,\bn_{n-1}),\\
   \prod_{i\in I_r}
  T_{\l_{i}}e^{h_n(\bl_n,\bn_{n-1})x_n}=&\,e^{2rx_n}\,\times\,e^{h_n(\bl_n,\bn_{n-1})x_n}
\end{split} \eeq}
On the other hand, operator $T_{\l_{i}}$ shifts the integration contour in \rf{p7}  in such a way that the conditions \rf{p1a} are satisfied with the replacement of $\l_i\equiv\g_{n,i}$ by $\l_i+2$, that is the shifted contour should separate set of poles:
 \beq \label{p9}\begin{split} \nu_{a}&=\l_{j}+g+2k,\qquad k=0,1,\ldots \qquad j\not\in I_r,
	\\ \nu_{a}&=\l_{j}+g+2k+2,\qquad k=0,1,\ldots \qquad j\in I_r
\end{split} \eeq
from
 \beq \label{p10}\begin{split} \nu_{a}&=\l_{j}-g-2k,\qquad k=0,1,\ldots \qquad j\not\in I_r,
	\\ \nu_{a}&=\l_{j}-g-2k+2,\qquad k=0,1,\ldots \qquad j\in I_r
\end{split} \eeq
Assume that all $\l_j$ are pure imaginary (or have a small real part, $\Re \l_j<(g-1)/2$).

Then the contour
\beq \label{p11}
C:\ \Re \nu_a=c, \qquad -g+2<c<g,\qquad a=1,\ldots, n-1
\eeq
satisfies the conditions \rf{p10}, \rf{p11} for any set $I_r$. Moreover, this contour may be obtained by such continuous deformation of imaginary plane $i\R^{n-1}$, which transforms the original configuration \rf{p1a} of poles into \rf{p9}--\rf{p10}. The picture below demonstrate the position of the poles in the integral $J_{I_r}$. Here $i\in I_r$ and
$j\not\in I_r$. The grey region describes the position of the contour $C$, specified for a variable $\nu_a$. This position  is common for all possible subsets $I_r\subset\nn$.
\bigskip

	\begin{tikzpicture}
\fill[gray!10] (0.55,-2.5) rectangle (1.45,2.5);
\fill(-3.5,1) circle(3pt); \node[below] at (-3.5,1)	
{$\lambda_j-g-2$};	
\fill(-1.5,1) circle(3pt); \node[below] at (-1.5,1)	
{$\lambda_j-g$};
\fill(-5.5,1) circle(3pt); \node[below] at (-5.5,1)	
{$\lambda_j-g-4$};
\fill(3.5,1) circle(3pt); \node[below] at (3.5,1)	
{$\lambda_j+g+2$};	
\fill(1.5,1) circle(3pt); \node[below] at (1.5,1)	
{$\lambda_j+g$};
\fill(5.5,1) circle(3pt); \node[below] at (5.5,1)	
{$\lambda_j+g+4$};
%\fill(7.5,1) circle(3pt); \node[below] at (7.5,1)	
%{$\lambda_i+g+6$};
%\fill(-7.5,3) circle(3pt); \node[below] at (-7.5,3)	
%{$\lambda_i-g-6$};

\fill(-3.5,-1) circle(3pt); \node[below] at (-3.5,-1)	
{$\lambda_i-g-2$};	
\fill(-1.5,-1) circle(3pt); \node[below] at (-1.5,-1)	
{$\lambda_i-g$};
%\fill(-5.5,-1) circle(3pt); \node[below] at (-5.5,-1)	
%{$\lambda_i-g-4$};
\fill(3.5,-1) circle(3pt); \node[below] at (3.5,-1)	
{$\lambda_i+g+2$};	
\fill(0.5,-1) circle(3pt); \node[below] at (0.5,-1)	
{$\lambda_i-g+2$};
\fill(5.5,-1) circle(3pt); \node[below] at (5.5,-1)	
{$\lambda_i+g+4$};
\fill(7.5,-1) circle(3pt); \node[below] at (7.5,-1)	
{$\lambda_j+g+6$};
%\fill(-7.5,-3) circle(3pt); \node[below] at (-7.5,-3)
%\draw[thin](-4,0)--(4,0);
\draw[dashed] (0.5,-2.5)--(0.5,2.5);
\draw[dashed] (1.5,-2.5)--(1.5,2.5);
\node[left] at (0.5,-2.5)	{$\operatorname{Re} \nu_a=-g+2$};
\node[right] at (1.5,-2.5)	{$\operatorname{Re} \nu_a=g$};
\draw[very thick,->] (1,-2.5)--(1,2.5);
\node[right] at (1,2.5) {$C$};
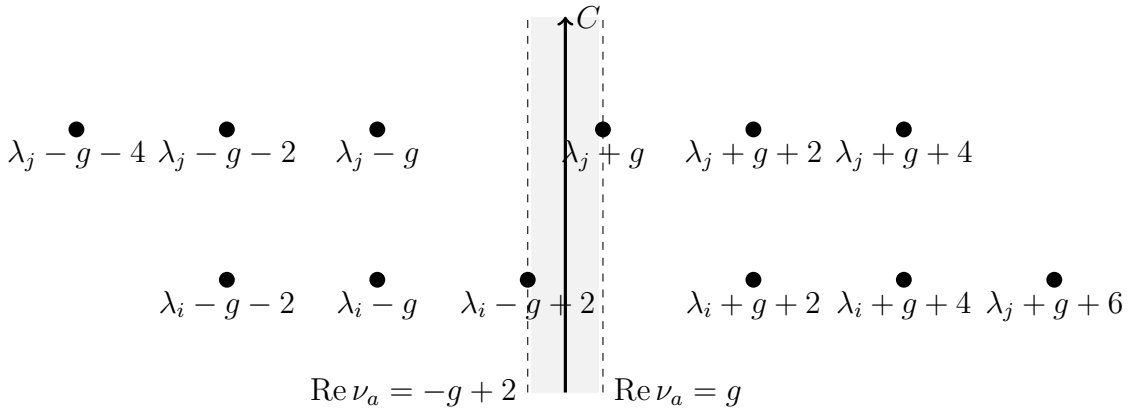
\node[below] at (0,-3.5) {Fig. 1 Contour $C$ for the variable $\nu_a$};
\end{tikzpicture}
\bigskip

Since the contour $C$ does not depend on a set $I_r$, we can permute integration and summation procedures  so that
 \beq \label{p12}\begin{split}
  &{\cal H}_r^{(g)}(\bl_n)\Phi_{\bl_n}^{(g)}(\x_n)=\sum_{I_r\subset\nn}J_{I_r}=\\&e^{2rx_n}\int\limits_{C}S_r(\bl_n,\bn_{n-1})
  \mu(\bn_{n-1})d\bn_{n-1}K(\bl_n,\bn_{n-1})e^{h_n(\bl)x_n}
  \Phi_{\bn_{n-1}}^{(g)}({\bm x}_{n-1}), 	
 	\end{split}
 \eeq
  where
  \beq\label{p13}
  S_r(\bl_n,\bn_{n-1})=\sum_{I_r\subset\nn}\prod_{i\in I_r}\left(\prod_{j\in \nn\setminus I_r}\frac{\l_i-\l_j+2-2g}{\l_i-\l_j}\prod_{a=1}^{N-1}\frac{\l_i-\nu_a+g}{\l_i-\nu_a+2-g}\right)
  \eeq
% Here by $\nn$ we denote the set $\nn=\{1,\ldots,n\}$ and by $I_r$ any subset of $\nn$ of cardinality $r$.
  Note that we eliminate the sign in \rf{p6} by rewriting the fractions.

 Define a similar sum,
 \beq\label{p14}
 \tS_r(\bn_{n-1},\bl_n,)=\sum_{A_r\subset\nnn}\prod_{a\in A_r}\left(\prod_{b\in \nnn\setminus A_r}
 \frac{\nu_a-\nu_b+2g-2}{\nu_a-\nu_b}\prod_{i=1}^{N}\frac{\l_i-\nu_a+g}{\l_i-\nu_a+2-g}\right)
 \eeq
 We need the following generalization of recurrence relations on binomial coefficients.
 \begin{lemma}\label{lemma1} For any $r=1,\ldots,n$
 	\beq\label{p15}
 S_r(\bl_n,\bn_{n-1})=\tS_{r-1}(\bn_{n-1},\bl_n,)+\tS_r(\bn_{n-1},\bl_n,)	
 	\eeq 	
 	\end{lemma}
 {\it Proof} is given  in Appendix.

\smallskip\noindent
 Using Lemma \ref{lemma1}, we rewrite \rf{p12} as
 \beq \label{p16}\begin{split}
 	&{\cal H}_r^{(g)}
 (\bl_n)\Phi_{\bl_n}^{(g)}(\x_n)=\\&e^{2rx_n}\int\limits_{C}\tS_{r-1}(\bn_{n-1},\bl_n,)
 	\mu(\bn_{n-1})d\bn_{n-1}K(\bl_n,\bn_{n-1})e^{h_n(\bl)x_n}
 	\Phi_{\bn_{n-1}}^{(g)}({\bm x}_{n-1})+\\
 	&e^{2rx_n}\int\limits_{C}\tS_r(\bn_{n-1},\bl_n,)
 	\mu(\bn_{n-1})d\bn_{n-1}K(\bl_n,\bn_{n-1})e^{h_n(\bl)x_n}
 	\Phi_{\bn_{n-1}}^{(g)}({\bm x}_{n-1}) 	
 \end{split}
 \eeq
  and apply to each occurring summand the same procedure in opposite direction.
  Namely, for any subset $A_r\subset \nnn$  of cardinality $r$ in the integral
  \beq\label{p16a}\begin{split}
  J'_{A_r}=\int\limits_{C}\prod_{a\in A_r}\left(\prod_{b\in \nnn\setminus A_r}
  \frac{\nu_a-\nu_b+2g-2}{\nu_a-\nu_b}\prod_{i=1}^{N}\frac{\l_i-\nu_a+g}{\l_i-\nu_a+2-g}\right)\\
 	\mu(\bn_{n-1})d\bn_{n-1}K(\bl_n,\bn_{n-1})e^{h_n(\bl)x_n}
 \Phi_{\bn_{n-1}}^{(g)}({\bm x}_{n-1})
  \end{split}\eeq
 we perform the change of integration variables
 \beq\label{p16b}
 \nu_a\to\nu_a+2,\qquad a\in A_r
 \eeq
  Functional relations on Gamma functions imply the relations
  \beq \label{p17}\begin{split}
  	 \prod_{a\in A_r}
  	T_{\nu_{a}}K(\bl_n,\bn_{n-1})=&\prod_{a\in A_r}\prod_{i=1}^{N}\frac{\nu_{a}-\l_i+g}{\l_i-\nu_{a}+g-2}K(\bl_n,\bn_{n-1}),\\
  	\prod_{a\in A_r}T_{\nu_{a}}\mu(\bn_{n-1})d\bn_{n-1}=&\prod_{a\in A_r, b\not\in A_r}\frac{\nu_a-\nu_b+2}{\nu_a-\nu_b}\cdot
  	\frac{\nu_a-\nu_b+2-2 g}{\nu_a-\nu_b+2g}\mu(\bn_{n-1})d\bn_{n-1}\,\\
  	\prod_{a\in A_r}
  	T_{\nu_{a}}e^{h_n(\bl_n,\bn_{n-1})x_n}=&e^{-2rx_n}\,\times\,e^{h_n(\bl_n,\bn_{n-1})x_n}
  \end{split} \eeq
We also have
\beq \label{p18}\begin{split}
  \prod_{a\in A_r}
 T_{\nu_{a}}\prod_{a\in A_r}\left(\prod_{b\in \nnn\setminus A_r}
 \frac{\nu_a-\nu_b+2g-2}{\nu_a-\nu_b}\prod_{i=1}^{N}\frac{\l_i-\nu_a+g}{\l_i-\nu_a+2-g}\right)=\\
 \prod_{a\in A_r}\left(\prod_{b\in \nnn\setminus A_r}
 \frac{\nu_a-\nu_b+2g}{\nu_a-\nu_b+2}\prod_{i=1}^{N}\frac{\l_i-\nu_a+g-2}{\l_i-\nu_a-g}\right).
 \end{split}\eeq
Combining  \rf{p17} and \rf{p18} we see that
\beq \label{p19}\begin{split}
J'_{A_r}=&(-1)^{rN}e^{2rx_N}\int\limits_{\tilde{C}}
\mu(\bn_{n-1})d\bn_{n-1}K(\bl_n,\bn_{n-1})e^{h_n(\bl)x_n}\\ &\prod_{a\in A_r\, b \in \nnn\setminus A_r}
\frac{\nu_a-\nu_b+2-2g}{\nu_a-\nu_b}
 \prod_{a\in A_r}
T_{\nu_{a}}\Phi_{\bn_{n-1}}^{(g)}({\bm x}_{n-1})
\end{split}\eeq
where the contour $\tilde{C}$ is the deformation of the contour $C$ according to the change of variables \rf{p16b}.
 In the assumption $g>1$ we may choose again $\tilde{C}=i\R^{n-1}$.
 Analogously, for any subset $A_{r-1}\subset \nnn$ the integral
  \beq\label{p20}\begin{split}
 	J'_{A_{r-1}}=\int\limits_{C}\prod_{a\in A_{r-1}}\left(\prod_{b\in \nnn\setminus A_{r-1}}
 	\frac{\nu_a-\nu_b+2g-2}{\nu_a-\nu_b}\prod_{i=1}^{N}\frac{\l_i-\nu_a+g}{\l_i-\nu_a+2-g}\right)\\
 	\mu(\bn_{n-1})d\bn_{n-1}K(\bl_n,\bn_{n-1})e^{h_n(\bl)x_n}
 	\Phi_{\bn_{n-1}}^{(g)}({\bm x}_{n-1})
 \end{split}\eeq
evaluates as
\beq \begin{split}\label{p21}
&J'_{A_{r-1}}=(-1)^{(r-1)N}e^{2(r-1)x_N}\int\limits_{i\R^{n-1}}
\mu(\bn_{n-1})d\bn_{n-1}K(\bl_n,\bn_{n-1})e^{h_n(\bl)x_n}\\
 &\prod_{a\in A_{r-1}\, b \in \nnn\setminus A_{r-1}}
\frac{\nu_a-\nu_b+2-2g}{\nu_a-\nu_b}
\prod_{a\in A_{r-1}}
T_{\nu_{a}}\Phi_{\bn_{n-1}}^{(g)}({\bm x}_{n-1})
\end{split}\eeq
Combining \rf{p15}, \rf{p20} and  \rf{p21} we arrive to the proof of Proposition \ref{proposition1}. \hfill{$\Box$}

 In this way we finish the proof of Theorem 1.
 \begin{remark} Note the difference on conditions for the coupling constant $g$ which guarantee the eigenfunction relations for Sutherland Hamiltonian $H_2$, see \rf{i1} and dual Hamiltonians $D_r^{(g)}$, see \rf{s2}. According to
 {\rm\cite{KK}}, the function \rf{i2}--\rf{i3} is an eigenfunction for the Sutherland Hamiltonian $H_2$ for any $g>0$. In this paper the  proof of Theorem 1 requires the condition $g>1$, which guarantees the existence of the common contour $C$ for all summands  $J_{I_r}$ in \rf{p6} due to the condition \rf{p11} which is nonempty only for $g>1$.
  	\end{remark}

\section{Example: Dual spectral problem for $n=2$}
From (\ref{i4}) in the case of $n=2$ one has two dual Hamiltonians
\beq\label{dh1}
\begin{split}
&{\cal H}_1^{(g)}(\l_1,\l_2)=-\frac{1}{\l_1-\l_2}
\Big\{(\l_1-\l_2+2-2g)e^{2\d_{\l_1}}+
(\l_1-\l_2-2+2g)e^{2\d_{\l_2}}\Big\},\\
&{\cal H}_2^{(g)}(\l_1,\l_2)=e^{2(\d_{\l_1}+\d_{\l_2})}
\end{split}
\eeq
In \cite{KK} we consider the simplest case of (\ref{i3}):
\beq\label{ac0}
\begin{split}
\Phi^{(g)}_{\l_1,\l_2}(x_1,x_2)&=\int\limits_{\imath\R}
\exp\Big\{(\l_1+\l_2-\g)x_2+\g x\Big\} \\ \times \, &
\Gamma\Big(\frac{\g-\l_1+g}{2}\Big)
\Gamma\Big(\frac{\l_1-\g+g}{2}\Big)
\Gamma\Big(\frac{\g-\l_2+g}{2}\Big)
\Gamma\Big(\frac{\l_2-\g+g}{2}\Big)
d\gamma,
\end{split}
\eeq
One can calculate integral (\ref{ac0}) in explicit terms.
Performing in (\ref{ac0}) the shift
$\g=s+\l_1$ of the integration variable, we rewrite the wave function as
\beq\label{ex0}
\Phi^{(g)}_{\l_1,\l_2}(x_1,x_2)=
\exp\Big\{\frac{1}{2}(\l_1+\l_2)(x_1+x_2)\Big\}\phi^{(g)}_{\l_1-\l_2}(x_1-x_2),
\eeq
where the $\frak{sl}(2,\R)$--part $\phi^{(g)}_{\l_1-\l_2}(x_1-x_2)$ has the form
\beq\label{ex1}
\phi^{(g)}_\l(x)=e^{\frac{\l x}{2}}\int\limits_{\imath\R}
\Gamma\Big(\frac{g+s}{2}\Big)
\Gamma\Big(\frac{g-s}{2}\Big)
\Gamma\Big(\frac{g+\l+s}{2}\Big)
\Gamma\Big(\frac{g-\l-s}{2}\Big)
\exp\{s x\}d s.
\eeq
Calculation of (\ref{ex1}) by residues yields up to numerical factor
\beq\label{mb4}
\phi^{(g)}_\l(x)=
\Gamma({\textstyle g-\frac{\l}{2}})
\Gamma({\textstyle g+\frac{\l}{2}})\cdot \sinh^{\frac{1}{2}-g} x\,P^{\frac{1}{2}-g}_{\frac{\l}{2}-\frac{1}{2}}(\cosh x),
\eeq
where $P^{\frac{1}{2}-g}_{\frac{\l}{2}-\frac{1}{2}}(\cosh x)$ is the Legendre function of the first kind.
Thus the wave function (\ref{ac0}) acquires the final form
\beq\label{ex6}
\begin{split}
\Phi^{(g)}_{\l_1,\l_2}(x_1,x_2)=
\sinh^{\frac{1}{2}-g}(x_1-x_2)\cdot
\exp\Big\{\frac{1}{2}(\l_1+\l_2)(x_1+x_2)\Big\}\\
\times\,
\Gamma({\textstyle g-\frac{\l_1-\l_2}{2}})
\Gamma({\textstyle g+\frac{\l_1-\l_2}{2}})%\cdot
 P^{\frac{1}{2}-g}_{\frac{\l_1-\l_2}{2}-\frac{1}{2}}(\cosh (x_1-x_2)).
\end{split}
\eeq
Since $P^{\frac{1}{2}-g}_{\frac{\l}{2}-\frac{1}{2}}(\cosh x)$ satisfies the relation
\beq
\frac{1}{\la}\Big\{(\la+2 g)e^{2\d_\la}+
(\la-2g)e^{-2\d_\l}\Big\}
P^{\frac{1}{2}-g}_{\frac{\la}{2}-\frac{1}{2}}(\cosh x)
=(e^x+e^{-x})P^{\frac{1}{2}-g}_{\frac{\l}{2}-\frac{1}{2}}(\cosh x)
\eeq
(see for example \cite[3.8(2)]{BE}), the function $\Phi^{(g)}_{\l_1,\l_2}(x_1,x_2)$ satisfies the spectral problem
\beq\label{aa1}
\begin{split}
&{\cal H}_1^{(g)}(\l_1,\l_2)\Phi^{(g)}_{\l_1,\l_2}(x_1,x_2)=(e^{2x_1}+e^{2x_2})
\Phi^{(g)}_{\l_1,\l_2}(x_1,x_2),\\
&{\cal H}_2^{(g)}(\l_1,\l_2)\Phi^{(g)}_{\l_1,\l_2}(x_1,x_2)=e^{2(x_1+x_2)}
\Phi^{(g)}_{\l_1,\l_2}(x_1,x_2),
\end{split}
\eeq
where dual operators ${\cal H}_1^{(g)},\,{\cal H}_2^{(g)}$ are given by (\ref{dh1}).

\appendix
\section{ Proof of Lemma \ref{lemma1}}
Rename first the variables. Put
\beq\label{a1}\begin{split} u_i&=\l_i+2-g, \qquad \ \ \, v_a=\nu_a, \qquad \a=2g-2,\\
\bu_n&=(u_1,\ldots,u_n),\qquad \bv_{n-1}=(v_1,\ldots, v_{n-1})	
\end{split}\eeq
 Then the identity  takes the form
\beq\label{a2}
S'_r(\bu_n,\bv_{n-1})=\tS'_{r-1}(\bv_{n-1},\bu_n,)+\tS'_r(\bv_{n-1},\bu_n,)	
\eeq
 where
\beq\begin{split}\label{a2a}
S'_r(\bu_n,\bv_{n-1})=\sum_{I_r\subset\nn}\prod_{i\in I_r}\left(\prod_{j\in \nn\setminus I_r}\frac{u_i-u_j-\a}{u_i-u_j}\prod_{a=1}^{N-1}\frac{u_i-v_a+\a}{u_i-v_a}\right),\\
\tS'_r(\bv_{n-1},\bv_n,)=\sum_{A_r\subset\nnn}\prod_{a\in A_r}\left(\prod_{b\in \nnn\setminus A_r}
\frac{v_a-v_b+\a}{v_a-v_b}\prod_{i=1}^{N}\frac{u_i-v_a+\a}{u_i-v_a}\right)
\end{split}\eeq
We prove \rf{a2} by induction on the number $n$ of variables $u_i$.
First we see that
\beq\label{a3}
S'_n(\bu_n,\bv_{n-1})=\tS'_{n-1}(\bv_{n-1},\bu_n)=\prod_{i\in \nn,\, a\in \nnn}\frac{u_i-v_a+\a}{u_i-v_a}
\eeq
 so that the equality \rf{a2} is valid for $n\leq r$. Next  assume that the equality \rf{a2} is valid for any number of variables $u_i$ less than $n$ (and any $r$) and prove it for the number of variables equal to $n$. The difference of LHS and RHS of \rf{a2} is symmetric with respect to permutations of $u_i$ and permutations of $v_a$. Thus the nominator of this rational function it is divisible by
  $$\Delta_{\bu_n}=\prod_{1\leq i<j\leq n}(u_i-u_j)\qquad \text{and}\qquad \Delta_{\bv_{n-1}}=\prod_{1\leq a<b\leq n-1}(v_a-v_b)$$
  and the only possible simple poles of this difference are
  $$ u_i=v_a.$$
  Consider  the difference of LHS and RHS of \rf{a2} as the function of $v_{n-1}$ and compute its residue at the point
   $$v_{n-1}=u_{n}=w.$$
   Denote by $E_n(w)$ the rational function
   \beq\label{a4}
   E_n(w)=\prod_{1\leq i<n}\frac{w-u_i-\a}{w-u_i}\prod_{1\leq a<n-1}\frac{w-v_a+\a}{w-v_a}
   \eeq
   One can observe the following equalities
   \beq\label{a5}\begin{split}
   \Res_{v_{n-1}=u_{n}=w}S'_r(\bu_n,\bv_{n-1})=\a E_n(w)S'_{r-1}(\bu_{n-1},\bv_{n-2}),\\
   \Res_{v_{n-1}=u_{n}=w}\tS'_r(\bv_{n-1},\bu_n)=\a E_n(w)\tS'_{r-1}(\bv_{n-2},\bu_{n-1})
   \end{split}\eeq
These equalities together with  induction assumptions imply that the residue of  the difference of LHS and RHS of \rf{a2}  at $v_{n-1}=u_{n}$ equals zero. By symmetries between $u_i$ and between $v_a$ all residues at $v_a-u_j$ also vanish.
 Thus this difference is a polynomial, symmetric with respect to $u_i$ and $v_j$ of degree zero, that is the constant.

 To compute this constant, tend all the variables $u_i$ and $v_a$ to infinity. It is equivalent to tend $\a$ to zero.
  But in this limit the equality \rf{a2} reduces to the recurrence relation on binomial coefficients,
   \beq \label{a6} \binom{n}{r}=\binom{n-1}{r-1}+\binom{n-1}{r} \eeq
   Thus the constant is zero. This proves the induction step. \hfill{$\Box$}

\section*{Funding} RFBR and NSFB according to the research project number 19-51-18006  (S.Kharchev); Russian Science Foundation, project No. 20-41-09009, used for the proof of the statements of Section  3  (S. Khoroshkin); section 2 was prepared within the framework of the HSE University Basic Research Program (S.Khor.)

	\end{document}